# Thèmes de parités dans l'évaluation aux entiers naturels de la fonction zêta et de cinq fonctions apparentées

par David Pouvreau[1]

**Résumé.** L'objet de cet article est le problème posé par l'évaluation aux entiers naturels de la fonction zêta et de cinq autres fonctions qui lui sont naturellement associées. Une approche relativement élémentaire est présentée, qui met ce problème encore partiellement ouvert en relation étroite avec quatre thèmes de parité : les notions de parité d'une fonction et de parité du degré d'un polynôme rejoignent ici les distinctions de parité concernant tant l'argument entier des six fonctions considérées que les entiers dont certaines puissances inverses sont sommées. La méthode adoptée a essentiellement pour but de familiariser avec ce problème les étudiants de mathématiques.

## 1. Introduction

On s'intéressera ici à la fonction zêta de Riemann, définie sur $]1; +\infty[$ par :

$$\zeta(x) = \sum_{n=1}^{+\infty} \frac{1}{n^x}$$

ainsi qu'à plusieurs autres fonctions apparentées. Ces dernières sont d'abord les deux fonctions $\alpha$ et $\beta$ directement liées à $\zeta$ par simples regroupements de termes, respectivement définies sur $]1; +\infty[$ par :

$$\alpha(x) = \sum_{m=1}^{+\infty} \frac{1}{(2m)^x} \quad ; \quad \beta(x) = \sum_{m=0}^{+\infty} \frac{1}{(2m+1)^x}$$

Les fonctions apparentées considérées sont aussi celles correspondant aux séries alternées qui sont naturellement associées aux trois séries précédentes, respectivement définies sur $]0; +\infty[$ par :

$$\xi(x) = \sum_{n=1}^{+\infty} \frac{(-1)^{n-1}}{n^x} \quad ; \quad \varphi(x) = \sum_{m=1}^{+\infty} \frac{(-1)^{m-1}}{(2m)^x} \quad ; \quad \psi(x) = \sum_{m=0}^{+\infty} \frac{(-1)^m}{(2m+1)^x}$$

Le problème considéré ici est le calcul des images des entiers naturels par ces six fonctions.

Il a suscité une vaste littérature depuis les travaux initiés à son sujet au XVIIIe siècle. La connaissance de $\psi(1)$, égal à $\pi/4$, est même très antérieure puisqu'on peut non sans justifications la considérer comme datant au plus tard du début du XVe siècle : on la trouve sous une forme ne faisant certes pas l'usage de l'infini actuel dans les travaux de l'astronome indien Mādhava de Sangamagrama[2]. Le mathématicien et philosophe allemand Gottfried W. Leibniz fut le premier Européen à l'établir, en 1673. Le problème fameux de la valeur de $\zeta(2)$, connu sous le nom de « problème de Bâle », avait été posé dès 1644 par le mathématicien italien Pietro Mengoli. Il fallut attendre 1735 pour le voir résolu par le mathématicien suisse Leonhard Euler : il trouva la valeur $\pi^2/6$ par une méthode ingénieuse mais logiquement incorrecte consistant à appliquer au développement en série entière de $\sin(x)/x$ un argument sur les relations entre racines et coefficients des polynômes, et détermina les valeurs de $\zeta(4)$, $\zeta(6)$, $\zeta(8)$, $\zeta(10)$ et $\zeta(12)$ par la même méthode[3]. Les progrès de l'analyse, surtout la théorie des séries de Fourier, ont permis au XIX$^e$ siècle de trouver les valeurs des cinq premières fonctions en les entiers pairs, et celles de la dernière aux entiers impairs. Mais les valeurs des cinq premières fonctions aux entiers impairs et celles de la sixième aux entiers pairs demeurent éminemment problématiques. Bien que plusieurs découvertes importantes à ce sujet aient été faites au cours des dernières décennies[4], le problème que posent ces valeurs demeure largement ouvert et fait toujours l'objet d'actives recherches.

---

[1] Professeur agrégé de mathématiques et docteur en histoire des sciences. Institut Alexander Grothendieck (Université de Montpellier) et Université de Mayotte. Email : david_pouvreau@orange.fr

[2] Les aspects historiques et mathématiques de cette découverte indienne sont discutés dans Pouvreau D., « Sur l'accélération de la convergence de la série de Mādhava-Leibniz », *Quadrature*, n° 97, 2015, pp. 17-25.

[3] Euler L., « De summis serierum reciprocarum », *Commentarii academiae, scientarum Petropolitanae* 7, 1740, pp. 123-134. Pour une traduction en anglais de cet article : https://arxiv.org/pdf/math/0506415.pdf

Une approche relativement élémentaire en est présentée ici, qui le met en relation étroite avec quatre thèmes de parité : les notions de parité d'une fonction et de parité du degré d'un polynôme rejoignent ici les distinctions de parité concernant tant l'argument entier des six fonctions considérées que les entiers dont certaines puissances inverses sont sommées. La méthode adoptée a essentiellement pour but de familiariser avec le problème considéré un public d'étudiants de mathématiques.

## 2. Prérequis : les deux théorèmes de Dirichlet sur les séries de Fourier

Tous les résultats de cet article découlent de deux propositions dont les démonstrations reposent sur l'utilisation de deux théorèmes qui seront ici dits « de Dirichlet » ; même si historiquement, seul le premier peut être attribué au mathématicien allemand Johann P.G. Lejeune Dirichlet, qui l'énonça et le démontra partiellement en 1829. Ils seront utilisés dans la version suivante :

Soit $F$ une fonction réelle continue par morceaux sur $\mathbb{R}$ et $T$-périodique. On pose $\omega = 2\pi/T$.
Pour tout $n \in \mathbb{N}$, on note :

$$a_n = \frac{2}{T}\int_0^T F(t)\cos(n\omega t)\,dt \quad ; \quad b_n = \frac{2}{T}\int_0^T F(t)\sin(n\omega t)\,dt$$

On considère la suite $(S_N)_{N\geq 1}$ dont les termes sont les fonctions définies sur $\mathbb{R}$ par :

$$S_N(t) = \frac{a_0}{2} + \sum_{n=1}^N (a_n\cos(n\omega t) + b_n\sin(n\omega t))$$

(1) si $F$ est de classe $C^1$ par morceaux sur $\mathbb{R}$, alors la suite $(S_N)_{N\geq 1}$ converge simplement sur $\mathbb{R}$ et sa limite, appelée « série de Fourier » de $F$, est la fonction $\tilde{F}$ définie sur $\mathbb{R}$ par :

$$\tilde{F}(t) = \frac{1}{2}\left(\lim_{\substack{x\to t\\x<t}} F(x) + \lim_{\substack{x\to t\\x>t}} F(x)\right)$$

(2) si $F$ est continue et de classe $C^1$ par morceaux sur $\mathbb{R}$, alors la suite $(S_N)_{N\geq 1}$ converge uniformément sur $\mathbb{R}$ et sa limite (série de Fourier de $F$) est la fonction $F$ elle-même.

Concernant les nombres $a_n$ et $b_n$ (« coefficients de Fourier » de $F$), on établit les résultats suivants :
si $F$ est paire, alors pour tout $n \in \mathbb{N}$ :

$$a_n = \frac{4}{T}\int_0^{\frac{T}{2}} F(t)\cos(n\omega t)\,dt \quad ; \quad b_n = 0$$

si $F$ est impaire, alors pour tout $n \in \mathbb{N}$ :

$$a_n = 0 \quad ; \quad b_n = \frac{4}{T}\int_0^{\frac{T}{2}} F(t)\sin(n\omega t)\,dt$$

---

[4] Fischler S., « Irrationalité des valeurs de zêta », Séminaire Bourbaki (2002-2003), vol. 45, pp. 27-62.

## 3. Images des entiers naturels pairs non nuls par $\zeta$, $\xi$ et $\varphi$

Dans l'approche du problème proposée ici, le calcul des images des entiers naturels pairs non nuls par $\zeta$, $\xi$ et $\varphi$ repose sur la proposition suivante :

**Proposition 1**

Soit $p \in \mathbb{N}^*$. Pour toute fonction réelle $f$ polynomiale et de degré inférieur ou égal à $2p$ sur $[0; \pi]$ :

$$f(0) = \frac{1}{\pi} \int_0^\pi f(t) \, dt - \frac{2}{\pi} \sum_{k=1}^p (-1)^{k-1} \left( f^{(2k-1)}(\pi) \, \xi(2k) + f^{(2k-1)}(0) \, \zeta(2k) \right)$$

Démonstration. Soit $F$ la fonction $2\pi$-périodique et paire qui coïncide avec $f$ sur $[0; \pi]$. $F$ est continue et de classe $C^1$ par morceaux. D'après le théorème de Dirichlet, la série de Fourier de $F$ converge donc uniformément vers $F$. On en déduit que pour tout $t \in [0; \pi]$ :

$$f(t) = \frac{1}{\pi} \int_0^\pi f(t) \, dt + \sum_{n=1}^{+\infty} a_n \cos(nt)$$

où :

$$\forall n \in \mathbb{N}^*, \ a_n = \frac{2}{\pi} \int_0^\pi f(t) \cos(nt) \, dt$$

Par conséquent :

$$f(0) = \frac{1}{\pi} \int_0^\pi f(t) \, dt + \sum_{n=1}^{+\infty} a_n$$

Or, pour tout $n \in \mathbb{N}^*$ fixé, on établit par récurrence sur $j \in \mathbb{N}^*$ que :

$$a_n = \frac{2}{\pi} \sum_{k=1}^j (-1)^{k-1} \frac{(-1)^n f^{(2k-1)}(\pi) - f^{(2k-1)}(0)}{n^{2k}} + \frac{2(-1)^j}{\pi n^{2j}} \int_0^\pi f^{(2j)}(t) \cos(nt) \, dt$$

En effet, deux intégrations par parties successives permettent d'obtenir :

$$a_n = -\frac{2}{\pi n} \int_0^\pi f'(t) \sin(nt) \, dt = \frac{2}{\pi n^2} \left( (-1)^n f'(\pi) - f'(0) \right) - \frac{2}{\pi n^2} \int_0^\pi f''(t) \cos(nt) \, dt$$

Ce qui montre que l'identité est satisfaite pour $j = 1$. Supposons que l'identité est satisfaite pour un certain $j \in \mathbb{N}^*$. On obtient en intégrant par parties deux fois que :

$$\int_0^\pi f^{(2j)}(t) \cos(nt) \, dt = -\frac{1}{n} \int_0^\pi f^{(2j+1)}(t) \sin(nt) \, dt$$
$$= \frac{1}{n^2} \left( (-1)^n f^{(2j+1)}(\pi) - f^{(2j+1)}(0) \right) - \frac{1}{n^2} \int_0^\pi f^{(2j+2)}(t) \cos(nt) \, dt$$

En utilisant l'hypothèse de récurrence, on en déduit que l'identité est satisfaite au rang $j + 1$. L'identité annoncée, vraie par récurrence pour tout $j \in \mathbb{N}^*$, est donc vraie en particulier en $j = p$ ; c'est-à-dire que :

$$a_n = \frac{2}{\pi} \sum_{k=1}^p (-1)^{k-1} \frac{(-1)^n f^{(2k-1)}(\pi) - f^{(2k-1)}(0)}{n^{2k}} + \frac{2(-1)^p}{\pi n^{2p}} \int_0^\pi f^{(2p)}(t) \cos(nt) \, dt$$

Or, $f$ étant polynomiale de degré inférieur ou égal à $2p$, la fonction $f^{(2p)}$ est constante sur $[0; \pi]$. On en déduit que :

$$\int_0^\pi f^{(2p)}(t)\cos(nt)\,dt = 0$$

D'où résulte :

$$a_n = \frac{2}{\pi}\sum_{k=1}^{p}(-1)^{k-1}\frac{(-1)^n f^{(2k-1)}(\pi) - f^{(2k-1)}(0)}{n^{2k}}$$

On obtient donc finalement :

$$f(0) = \frac{1}{\pi}\int_0^\pi f(t)\,dt + \sum_{n=1}^{+\infty}\frac{2}{\pi}\sum_{k=1}^{p}(-1)^{k-1}\frac{(-1)^n f^{(2k-1)}(\pi) - f^{(2k-1)}(0)}{n^{2k}}$$

$$= \frac{1}{\pi}\int_0^\pi f(t)\,dt - \frac{2}{\pi}\sum_{k=1}^{p}(-1)^{k-1}\sum_{n=1}^{+\infty}\left(\frac{(-1)^{n-1}}{n^{2k}}f^{(2k-1)}(\pi) + \frac{1}{n^{2k}}f^{(2k-1)}(0)\right)$$

Et cette identité n'est autre que celle annoncée dans la proposition 1.

Pour obtenir la valeur de $\xi(2p)$ et $\zeta(2p)$ pour tout $p \in \mathbb{N}^*$, il suffit alors de choisir convenablement la fonction $f$ considérée dans la proposition 1.

**Proposition 2**

Il existe une suite $(A_p) \in \mathbb{Q}^{\mathbb{N}^*}$ telle que $\xi(2p) = A_p \pi^{2p}$ pour tout $p \in \mathbb{N}^*$.

Les termes de cette suite sont déterminés par $A_1 = 1/12$ et par la relation de récurrence multiple :

$$\forall p \geq 2,\ A_p = \frac{(-1)^{p-1}}{(2p)!}\left(\frac{1}{2(2p+1)} - \sum_{k=1}^{p-1}(-1)^{k-1}\frac{(2p)!}{(2p+1-2k)!}A_k\right)$$

<u>Démonstration</u>. Soit $p \in \mathbb{N}^*$. On applique la proposition 1 avec $f : t \mapsto t^{2p}$. On obtient alors :

$$0 = \frac{1}{\pi}\frac{\pi^{2p+1}}{2p+1} - \frac{2}{\pi}\sum_{k=1}^{p}(-1)^{k-1}\frac{(2p)!}{(2p+1-2k)!}\pi^{(2p+1)-2k}\xi(2k)$$

Soit aussi :

$$\sum_{k=1}^{p}(-1)^{k-1}\frac{(2p)!}{(2p+1-2k)!}\frac{1}{\pi^{2k}}\xi(2k) = \frac{1}{2(2p+1)}$$

En posant $A_q = \xi(2q)/\pi^{2q}$ pour tout $q \in \mathbb{N}^*$, on obtient la relation exprimant $A_p$ lorsque $p \geq 2$.

Dans le cas où $p = 1$, on obtient de plus avec la relation précédemment obtenue que $\xi(2) = \pi^2/12$, donc que $\xi(2)/\pi^2 = A_1 = 1/12 \in \mathbb{Q}$.

Supposons alors que, pour un $p \in \mathbb{N}^*$ fixé, $A_j \in \mathbb{Q}$ pour tout $j \in [\![1;p]\!]$. Il est clair que

$$A_{p+1} = \frac{\xi(2p+2)}{\pi^{2p+2}} = \frac{(-1)^p}{(2p+2)!}\left(\frac{1}{2(2p+3)} - \sum_{k=1}^{p}(-1)^{k-1}\frac{(2p+2)!}{(2p+3-2k)!}A_k\right) \in \mathbb{Q}$$

Par récurrence forte sur $p$, on en déduit que $A_p \in \mathbb{Q}$ pour tout $p \in \mathbb{N}^*$.

**Proposition 3**

Pour tout $p \in \mathbb{N}^*$ :

$$\varphi(2p) = \frac{A_p}{4^p}\pi^{2p}$$

<u>Démonstration</u> : Il suffit d'observer que :

$$\varphi(2p) = \sum_{m=1}^{+\infty}\frac{(-1)^{m-1}}{(2m)^{2p}} = \frac{1}{4^p}\sum_{m=1}^{+\infty}\frac{(-1)^{m-1}}{m^{2p}} = \frac{1}{4^p}\xi(2p)$$

**Proposition 4**

Il existe une suite $(B_p) \in \mathbb{Q}^{\mathbb{N}^*}$ telle que $\zeta(2p) = B_p\pi^{2p}$ pour tout $p \in \mathbb{N}^*$.

Les termes de cette suite sont déterminés par $B_1 = 1/6$ et par la relation de récurrence multiple :

$$\forall\, p \geq 2,\ B_p = \frac{(-1)^{p-1}}{(2p)!}\left(\frac{p}{2p+1} - \sum_{k=1}^{p-1}(-1)^{k-1}\frac{(2p)!}{(2p+1-2k)!}B_k\right)$$

On a de plus la relation :

$$\forall\, p \geq 1,\ B_p = \frac{4^p}{4^p - 2}A_p$$

<u>Démonstration</u>. Soit $p \in \mathbb{N}^*$. On applique encore la proposition 1, mais avec cette fois-ci la fonction $f$ définie sur $[0;\pi]$ par : $f(t) = (t-\pi)^{2p}$. On obtient alors :

$$\pi^{2p} = \frac{1}{\pi}\frac{\pi^{2p+1}}{2p+1} - \frac{2}{\pi}\sum_{k=1}^{p}(-1)^{k-1}\frac{(2p)!}{(2p+1-2k)!}(-1)^{2p-2k+1}\pi^{2p+1-2k}\zeta(2k)$$

Soit aussi :

$$\frac{p}{2p+1} = \sum_{k=1}^{p}(-1)^{k-1}\frac{(2p)!}{(2p+1-2k)!}\frac{1}{\pi^{2k}}\zeta(2k)$$

En posant $B_q = \zeta(2q)/\pi^{2q}$ pour tout $q \in \mathbb{N}^*$, on obtient la relation exprimant $B_p$ lorsque $p \geq 2$.

Dans le cas où $p = 1$, on obtient de plus avec la relation précédemment obtenue que $\zeta(2) = \pi^2/6$, donc que $\zeta(2)/\pi^2 = B_1 = 1/6 \in \mathbb{Q}$.

Supposons alors que, pour un $p \in \mathbb{N}^*$ fixé, $B_j \in \mathbb{Q}$ pour tout $j \in [\![1;p]\!]$. Il est clair que

$$B_{p+1} = \frac{\zeta(2p+2)}{\pi^{2p+2}} = \frac{(-1)^p}{(2p+2)!}\left(\frac{p+1}{2p+3} - \sum_{k=1}^{p}(-1)^{k-1}\frac{(2p+2)!}{(2p+3-2k)!}B_k\right) \in \mathbb{Q}$$

Par récurrence forte sur $p$, on en déduit que $B_p \in \mathbb{Q}$ pour tout $p \in \mathbb{N}^*$.

La dernière relation annoncée, qui lie $B_p$ à $A_p$, sera établie plus loin (proposition 5).

## 4. Images des entiers naturels non nuls par $\alpha$ et $\beta$

Le calcul des images des entiers naturels non nuls par les fonctions $\alpha$ et $\beta$ est élémentaire à partir de celui effectué avec les fonctions $\zeta$ et $\xi$. Leurs valeurs aux entiers pairs non nuls sont donc connues à partir des propositions 2 et 4 :

**Proposition 5**

Pour tout $p \in \mathbb{N}^*$ :

$$\alpha(p) = \frac{1}{2}(\zeta(p) - \xi(p)) \quad \beta(p) = \frac{1}{2}(\zeta(p) + \xi(p))$$

En particulier, pour tout $p \in \mathbb{N}^*$ :

$$\alpha(2p) = \frac{1}{4^p} B_p \pi^{2p} \quad ; \quad \beta(2p) = \frac{4^p - 1}{4^p} B_p \pi^{2p}$$

<u>Démonstration.</u> Pour tout $p \in \mathbb{N}^*$ :

$$\zeta(p) = \sum_{n=1}^{+\infty} \frac{1}{n^p} = \sum_{m=1}^{+\infty} \frac{1}{(2m)^p} + \sum_{m=0}^{+\infty} \frac{1}{(2m+1)^p} = \alpha(p) + \beta(p)$$

et

$$\xi(p) = \sum_{n=1}^{+\infty} \frac{(-1)^{n-1}}{n^p} = \sum_{m=1}^{+\infty} \frac{(-1)^{2m-1}}{(2m)^p} + \sum_{m=0}^{+\infty} \frac{(-1)^{2m}}{(2m+1)^p} = -\alpha(p) + \beta(p)$$

D'où les deux identités générales annoncées, puis leurs applications aux entiers pairs. Pour ces dernières, on utilise aussi les observations que :

$$\alpha(2p) = \sum_{m=1}^{+\infty} \frac{1}{(2m)^{2p}} = \frac{1}{4^p} \zeta(2p)$$

$$\zeta(2p) = \sum_{m=1}^{+\infty} \frac{1}{(2m)^{2p}} + \sum_{m=0}^{+\infty} \frac{1}{(2m+1)^{2p}} = \frac{1}{4^p} \zeta(2p) + \beta(2p)$$

dont résulte :

$$\beta(2p) = \frac{4^p - 1}{4^p} \zeta(2p)$$

C'est ce résultat qui justifie l'égalité annoncée à la fin de la proposition 3. En effet, on obtient par identification des deux expressions connues de $\beta(2p)$ la nouvelle identité :

$$\frac{1}{2}(B_p + A_p) = \frac{4^p - 1}{4^p} B_p$$

dont résulte :

$$A_p = \frac{4^p - 2}{4^p} B_p$$

## 5. Images des entiers naturels impairs par $\psi$

Une approche similaire à celle employée au 3. permet le calcul des images des entiers naturels impairs par $\psi$ :

**Proposition 6**

Soit $p \in \mathbb{N}$. Pour toute fonction réelle $f$ impaire, définie sur $[-\pi; \pi]$ et polynomiale de degré inférieur ou égal à $2p + 1$ sur cet intervalle :

$$f\left(\frac{\pi}{2}\right) = \frac{2}{\pi} \sum_{k=0}^{p} (-1)^k f^{(2k)}(\pi) \psi(2k+1)$$

<u>Démonstration</u>. Soit $F$ la fonction $2\pi$-périodique qui coïncide avec $f$ sur $]-\pi;\pi[$ et est telle que $F(n\pi) = 0$ pour tout $n \in \mathbb{N}$. $F$ est impaire et de classe $C^1$ par morceaux. La série de Fourier de $F$ converge donc simplement vers la fonction régularisée de $F$ d'après le théorème de Dirichlet. On en déduit que pour tout $t \in ]-\pi;\pi[$ :

$$f(t) = \sum_{n=1}^{+\infty} b_n \sin(nt)$$

où :

$$\forall\, n \in \mathbb{N}^*,\ b_n = \frac{2}{\pi} \int_0^\pi f(t) \sin(nt)\, dt$$

Par conséquent :

$$f\left(\frac{\pi}{2}\right) = \sum_{m=0}^{+\infty} (-1)^m b_{2m+1}$$

Or, pour tout $n \in \mathbb{N}^*$ fixé, on établit par récurrence sur $j \in \mathbb{N}$ que :

$$b_n = \frac{2}{\pi} \sum_{k=0}^{j} (-1)^{k-1} \frac{(-1)^n f^{(2k)}(\pi)}{n^{2k+1}} + \frac{2(-1)^j}{\pi n^{2j+1}} \int_0^\pi f^{(2j+1)}(t) \cos(nt)\, dt$$

En effet, on obtient en intégrant trois fois par parties et en tenant compte du fait que l'imparité de $f$ implique celle de toutes ses dérivées d'ordre pair (et donc leur nullité en 0) :

$$b_n = -\frac{2}{\pi n}(-1)^n f(\pi) + \frac{2}{\pi n} \int_0^\pi f'(t) \cos(nt)\, dt$$

$$= -\frac{2}{\pi n}(-1)^n f(\pi) + \frac{2}{\pi n^3}(-1)^n f''(\pi) - \frac{2}{\pi n^3} \int_0^\pi f^{(3)}(t) \cos(nt)\, dt$$

Ce qui montre que l'identité est satisfaite successivement pour $j = 0$ et $j = 1$. Supposons que l'identité annoncée vaut pour un certain $j \in \mathbb{N}$. On obtient, en intégrant deux fois par parties :

$$\int_0^\pi f^{(2j+1)}(t) \cos(nt)\, dt = -\frac{1}{n} \int_0^\pi f^{(2j+2)}(t) \sin(nt)\, dt$$

$$= \frac{1}{n^2}(-1)^n f^{(2j+2)}(\pi) - \frac{1}{n^2} \int_0^\pi f^{(2j+3)}(t) \cos(nt)\, dt$$

En utilisant l'hypothèse de récurrence, on en déduit que l'identité est satisfaite au rang $j + 1$. L'identité annoncée, vraie par récurrence pour tout $j \in \mathbb{N}$, est donc vraie en particulier en $j = p$ ; c'est-à-dire que :

$$b_n = \frac{2}{\pi} \sum_{k=0}^{p} (-1)^{k-1} \frac{(-1)^n f^{(2k)}(\pi)}{n^{2k+1}} + \frac{2(-1)^p}{\pi n^{2p+1}} \int_0^\pi f^{(2p+1)}(t) \cos(nt)\, dt$$

Or, $f$ étant polynomiale de degré inférieur ou égal à $2p + 1$ sur $[0;\pi]$, la fonction $f^{(2p+1)}$ y est constante. On en déduit :

$$\int_0^\pi f^{(2p+1)}(t) \cos(nt)\, dt = 0$$

D'où résulte :

$$b_n = \frac{2}{\pi} \sum_{k=0}^{p} (-1)^{k-1} \frac{(-1)^n f^{(2k)}(\pi)}{n^{2k+1}}$$

On obtient donc finalement :

$$f\left(\frac{\pi}{2}\right) = \sum_{m=0}^{+\infty} \frac{2(-1)^m}{\pi} \sum_{k=0}^{p} (-1)^{k-1} \frac{(-1)^{2m+1} f^{(2k)}(\pi)}{(2m+1)^{2k+1}}$$

$$= \frac{2}{\pi} \sum_{k=0}^{p} (-1)^k \left( \sum_{m=0}^{+\infty} \frac{(-1)^m}{(2m+1)^{2k+1}} f^{(2k)}(\pi) \right)$$

Et cette identité n'est autre que celle annoncée dans la proposition 6.

**Proposition 7**

Il existe une suite $(C_p) \in \mathbb{Q}^{\mathbb{N}}$ telle que $\psi(2p+1) = C_p \, \pi^{2p+1}$ pour tout $p \in \mathbb{N}$.

Les termes de cette suite sont déterminés par $C_0 = 1/4$ et par la relation de récurrence multiple :

$$\forall \, p \geq 1, \quad C_p = \frac{(-1)^p}{(2p+1)!} \left( \frac{1}{4^{p+1}} - \sum_{k=0}^{p-1} (-1)^k \frac{(2p+1)!}{(2p+1-2k)!} C_k \right)$$

<u>Démonstration</u>. Soit $p \in \mathbb{N}$. On applique la proposition 6 avec $f : t \mapsto t^{2p+1}$. On obtient alors :

$$\frac{\pi^{2p+1}}{2^{2p+1}} = \frac{2}{\pi} \sum_{k=0}^{p} (-1)^k \frac{(2p+1)!}{(2p+1-2k)!} \pi^{(2p+1)-2k} \psi(2k+1)$$

Soit aussi :

$$\sum_{k=0}^{p} (-1)^k \frac{(2p+1)!}{(2p+1-2k)!} \frac{1}{\pi^{2k+1}} \psi(2k+1) = \frac{1}{4^{p+1}}$$

En posant $C_q = \psi(2q+1)/\pi^{2q+1}$ pour tout $q \in \mathbb{N}$, on obtient la relation exprimant $C_p$ lorsque $p \geq 1$.

Dans le cas où $p = 0$, on obtient de plus avec la relation précédemment obtenue que $\psi(1) = \pi/4$, donc que $\psi(1)/\pi = C_0 = 1/4 \in \mathbb{Q}$.

Supposons alors que, pour un $p \in \mathbb{N}$ fixé, $C_j \in \mathbb{Q}$ pour tout $j \in [\![0; p]\!]$. Il est clair que

$$C_{p+1} = \frac{\psi(2p+3)}{\pi^{2p+3}} = \frac{(-1)^{p+1}}{(2p+3)!} \left( \frac{1}{4^{p+2}} - \sum_{k=0}^{p} (-1)^k \frac{(2p+3)!}{(2p+3-2k)!} C_k \right) \in \mathbb{Q}$$

Par récurrence forte sur $p$, on en déduit que $C_p \in \mathbb{Q}$ pour tout $p \in \mathbb{N}$.

## 6. Typologie liant la résolution du problème à quatre « niveaux » de parité

Tous les calculs possibles permettant de résoudre les problèmes initialement posés *au moyen de l'approche utilisée* ont été effectués. Concernant les fonctions envisagées au 3. et au 5., leur évaluation en d'autres réels de $[0; \pi[$ que 0 et $\pi/2$ ne permet pas de progresser – le lecteur pourra notamment vérifier que l'application à $\pi/2$ d'une fonction $f$ satisfaisant les hypothèses de la proposition 1 ne fournit aucun résultat nouveau. On retrouve ici exactement les cas problématiques actuellement encore non réglés, qui restent soumis aux investigations contemporaines ; c'est-à-dire celui du calcul des images des entiers impairs par $\zeta, \alpha, \beta, \xi$ et $\varphi$, et du calcul des images des entiers pairs par $\psi$. On ignore en fait toujours si, pour $p \in \mathbb{N}^*$ donné, le rapport à $\pi^{2p+1}$ des images de $2p+1$ par les cinq premières fonctions est un nombre rationnel, et la même ignorance vaut pour le rapport à $\pi^{2p}$ de $\psi(2p)$, pour tout $p \in \mathbb{N}$. En particulier, on ignore même si, hormis le cas de $\zeta(3)$, démontré irrationnel en 1978 par le mathématicien français Roger Apery, l'un de ces nombres arbitrairement choisi est bien irrationnel.

Le cas des images de 1 par les fonctions $\zeta, \alpha, \beta, \xi$ et $\varphi$ mérite enfin d'être relevé car sa spécificité interroge au regard de tout ce qui précède. Il est clair que $\zeta(1)$, $\alpha(1)$ et $\beta(1)$ n'existent pas. Il n'en va par contre pas de même de $\xi(1)$ et de $\varphi(1)$, mais il se trouve que leurs valeurs ne sont pas, contrairement à la logique qui se dégage des cas résolus, un multiple rationnel de $\pi$. En effet :

$$\xi(1) = \ln(2) \quad \text{et} \quad \varphi(1) = \frac{1}{2}\ln(2)$$

Il est d'ailleurs remarquable que la détermination de ces deux valeurs ne s'obtient pas par le même type de procédure fondée sur les séries de Fourier. La seconde se déduit sans difficulté de la première. La valeur de $\xi(1)$ se déduit du développement en série entière de $\ln(1 + x)$ et de l'argument de la convergence de la série en 1 par le critère spécial de convergence des séries alternées. Elle peut aussi être établie de manière moins sophistiquée comme suit :

$$\forall\, t \in \mathbb{R}, \forall\, n \in \mathbb{N}, \quad \sum_{k=0}^{n}(-1)^k t^k = \frac{1}{1+t} + (-1)^n \frac{t^{n+1}}{1+t}$$

D'où résulte par intégration sur $[0; 1]$ :

$$\sum_{k=0}^{n}\frac{(-1)^k}{k+1} = \ln(2) + (-1)^n \int_0^1 \frac{t^{n+1}}{1+t}\, dt$$

La valeur de $\xi(1)$ s'en déduit par passage à la limite lorsque $n$ tend vers $+\infty$, compte tenu de :

$$\forall\, n \in \mathbb{N}, \quad 0 \leq \int_0^1 \frac{t^{n+1}}{1+t}\, dt \leq \int_0^1 t^{n+1} dt = \frac{1}{n+2}$$

qui implique

$$\lim_{n \to +\infty} \int_0^1 \frac{t^{n+1}}{1+t}\, dt = 0$$

L'approche utilisée dans cet article a le mérite de révéler de manière relativement élémentaire deux aspects intimement liés auxquels se rapporte la résolution du problème considéré : (1) le déphasage des fonctions sinus et cosinus, dont la conséquence est manifeste dans l'inégalité du nombre de cas résolus par les propositions 1 et 6 ; (2) les liens profonds existant entre les résolutions cherchées et quatre situations différentes de parité : (a) la parité des fonctions $F$ utilisées ; (b) la parité du degré des restrictions polynomiales $f$ de ces fonctions à $[0; \pi[$ ; (c) la parité des entiers formant les arguments des six fonctions sur lesquels porte la question initialement posée ; (d) la parité des entiers dont les puissances inverses interviennent dans les sommes définissant les quatre fonctions $\alpha, \beta, \varphi$ et $\psi$.

Les deux tableaux suivants synthétisent une typologie des cas de résolution et des cas de non résolution relativement à ces quatre « niveaux » de parité :

| Arguments entiers naturels *pairs* comme puissances d'inverses d'entiers naturels *pairs* : | Arguments entiers naturels *impairs* comme puissances d'inverses d'entiers naturels *pairs* |
|---|---|
| Cas de $\alpha$ (série non alternée) résolu | Cas de $\alpha$ (série non alternée) non résolu |
| Cas de $\varphi$ (série alternée) résolu | Cas de $\varphi$ (série alternée) non résolu sauf $\varphi(1)$ |
| Arguments naturels *pairs* comme puissances d'inverses d'entiers naturels *impairs* | Arguments entiers naturels *impairs* comme puissances d'inverses d'entiers naturels *impairs* |
| Cas de $\beta$ (série non alternée) résolu | Cas de $\beta$ (série non alternée) non résolu |
| Cas de $\psi$ (série alternée) non résolu | Cas de $\psi$ (série alternée) résolu |

| Arguments entiers naturels *pairs* $2p \geq 2$ | Arguments entiers naturels *impairs* $2p + 1 \geq 1$ |
|---|---|
| Cas de $\zeta, \alpha, \beta, \xi$ et $\varphi$ résolus avec des fonctions périodiques *paires*, dont la restriction à $[0; \pi[$ est polynomiale et de degré *pair* | Cas de $\zeta, \alpha, \beta, \xi$ et $\varphi$ non résolus si $p \geq 1$ Cas particuliers si $p = 0$ : $\zeta(1), \alpha(1)$ et $\beta(1)$ non définis, $\varphi(1) = \frac{1}{2}\ln(2)$ et $\xi(1) = \ln(2)$ |
| Cas de $\psi$ non résolu | Cas de $\psi$ résolu avec des fonctions périodiques *impaires*, dont la restriction à $[0; \pi[$ est polynomiale et de degré *impair* |